\documentstyle{amsppt}
\tolerance 3000 \pagewidth{5.5in} \vsize7.0in
\magnification=\magstep1 \NoBlackBoxes \NoRunningHeads
\widestnumber \key{AAAAAAAAAAA} \topmatter
\author A. Iosevich and I. {\L}aba
\endauthor
\title $K$-distance sets, Falconer conjecture, and discrete
analogs
\endtitle
\date March 18, 2003
\enddate
\thanks Research of A. Iosevich supported in part by the NSF
Grant DMS00-87339 \endthanks
\thanks Research of I. {\L}aba supported in part by the NSERC
Grant 22R80520
\endthanks
\thanks 2000 Mathematics Subject Classification: 42B99
\endthanks
\address A. Iosevich, University of Missouri
\ email: iosevich \@ wolff.math.missouri.edu
\endaddress
\address I. {\L}aba, University of British Columbia
\ email: ilaba \@math.ubc.ca
\endaddress
\abstract In this paper we prove a series of results on the size
of distance sets corresponding to sets in the Euclidean space.
These distances are generated by bounded convex sets and the
results depend explicitly on the geometry of these sets. We also
use a diophantine mechanism to convert continuous results into
distance set estimates for discrete point sets. \endabstract
\endtopmatter
\document

\head Introduction \endhead

\vskip.125in

Let $E$ be a compact subset of ${\Bbb R}^d$. Let
$\Delta(E)=\{|x-y|:x,y \in E\}$, where $|\cdot|$ denotes the usual
Euclidean metric. The celebrated Falconer conjecture says that if
the Hausdorff dimension of $E$ is greater than $\frac{d}{2}$, then
$\Delta(E)$ has positive Lebesgue measure. Falconer
(\cite{Falconer86}) obtained this conclusion if the Hausdorff
dimension of $E$ is greater than $\frac{d+1}{2}$. His result was
improved by Bourgain in \cite{Bourgain94}. The best known result
in the plane is due to Tom Wolff who proved in \cite{Wolff99} that
the distance set has positive Lebesgue measure if the Hausdorff
dimension of $E$ is greater than $\frac{4}{3}$.

All these results are based on the curvature of the unit circle of
the Euclidean metric. For if the Euclidean metric is replaced by
the $l^{\infty}$ metric, for example, the situation becomes very
different. To see this, let $E=C_{2m} \times C_{2m}$, $m>1$ an
even integer, where $C_{2m}$ is the Cantor-type subset of $[0,1]$
consisting of numbers whose base $2m$ expansions contain only even
numbers. One can check that the distance set with respect to the
$l^{\infty}$ metric has measure $0$ for any $m$, whereas the
Hausdorff dimension of this set is $2\frac{\log(m)}{\log(2m)} \to
2$ as $m \to \infty$.

The curvature is not the end of the story. The aforementioned
results of Falconer, Bourgain and Wolff also used the smoothness
of the unit circle. Falconer did so explicitly by using the
formula for the Fourier transform of the characteristic function
of an Euclidean annulus, whereas Bourgain and Wolff utilized it
implicitly by using a reduction to circular averages which again
relies on asymptotics of the Fourier transform of the Lebesgue
measure on the unit circle of the distance which does not hold in
the absence of smoothness. See \cite{Mattila87} and
\cite{Sj\"olin93} for a background on these reductions.

The purpose of this paper to study the Falconer conjecture in the
absence of smoothness using geometric features of the Fourier
transforms of characteristic functions of convex sets in Euclidean
space. We also develop a conversion mechanism based on diophantine
approximation which allows us to obtain geometric combinatorial
results from analogous facts in a continuous setting. 

\head Continuous results: variants of the Falconer distance
problem \endhead

\vskip.125in

\definition{Definition}
Let $K$ be a bounded convex set in ${\Bbb R}^d$ symmetric about
the origin, and let ${||\cdot||}_K$ be the norm induced by $K$.
The $K$-distance set of a set $E\subset{\Bbb R}^d$ is the set
$\Delta_K(E)=\{{||x-y||}_K: x,y \in E\}$.
\enddefinition

Our main results are the following.

\proclaim{Theorem 0.1} Let $E \subset {[0,1]}^d$. Let $K$ be a
bounded convex set in ${\Bbb R}^d$ and let $\sigma_K$ denote the
Lebesgue measure on $\partial K$.

i) Suppose that $|\widehat{\sigma}_K(\xi)| \leq C
{|\xi|}^{-\gamma}$ for some $\gamma>0$, and the Hausdorff
dimension of $E$ is greater than $d-\gamma$. Then $\Delta_K(E)$
has positive Lebesgue measure.

ii) Suppose that $\int_{S^{d-1}} |\widehat{\sigma}_K(R \omega)|
d\omega \leq C{R}^{-\gamma}$, for some $\gamma>0$, and the
Hausdorff dimension of $E$ is greater than $d-\gamma$. Then
$\Delta_K(\rho E)$ has positive Lebesgue measure for almost every
$\rho \in S^{d-1}$, viewed (in the obvious way) as an element of
$SO(d)$ equipped with the natural measure.

iii) Suppose that
$$|\widehat{\sigma}_K(\xi)| \leq C \gamma({|\xi|}^{-1}), \tag0.1$$
where $\gamma$ is a convex function with $\gamma(0)=0$, and
$\sigma_K$ is the Lebesgue measure on $\partial K$. Let $E \subset
{[0,1]}^d$. Suppose that there exists a Borel measure on $E$ such
that
$$\int {|\widehat{\mu}(\xi)|}^2 \gamma({|\xi|}^{-1})d\xi<\infty.
\tag0.2$$

Then $\Delta_K(E)$ has positive Lebesgue measure.
\endproclaim

\proclaim{Corollary 0.2} Let $K$ be a bounded symmetric convex set
in ${\Bbb R}^d$.

i) Let $d=2$. Let $S_{\theta}=\sup_{x \in K} x \cdot \omega$,
where $\omega= (\cos(\theta), \sin(\theta))$, and denote by
$l(\theta,\epsilon)$ the length of the chord
$C(\theta,\varepsilon)=\left\{x \in K: x \cdot \omega= S_{\theta}-
\varepsilon \right\}$. Suppose that $\partial K$ has everywhere
non-vanishing curvature in the sense that there exists a positive
uniform constant $c$ such that
$$ l(\theta, \epsilon) \leq c \sqrt{\epsilon}. \tag0.3$$

Let $E \subset {[0,1]}^2$ be a set of Hausdorff dimension
$\alpha>\frac{1}{2}$. Then the Hausdorff dimension of
$\Delta_K(E)$ is at least $\min(1,\alpha-\frac{1}{2})$.  If
furthermore $\alpha>\frac{3}{2}$, then $\Delta_K(E)$ has positive
Lebesgue measure.

ii) Suppose that $E \subset {[0,1]}^d$ is a set of Hausdorff
dimension $\alpha>\frac{d+1}{2}$, and that $\partial K$ is smooth
and has non-vanishing Gaussian curvature.  Then $\Delta_K(E)$ has
positive Lebesgue measure.
\endproclaim

\proclaim{Corollary 0.3} Let $E \subset {[0,1]}^d$.

i) Suppose that $K$ is a symmetric convex polyhedron and that the
Hausdorff dimension of $E$ is greater than $1$. Then
$\Delta_K(\rho E)$ has positive Lebesgue measure for almost every
$\rho \in SO(d)$. Moreover, this result is sharp in the sense that
for every $\alpha<1$ there exists a set of Hausdorff dimension
$\alpha$ such that $\Delta_K(E)$ has Lebesgue measure $0$ with
respect to every convex body $K$.

ii) Suppose that $K$ is any bounded symmetric convex body in
${\Bbb R}^d$ and the Hausdorff dimension of $E$ is greater than
$\frac{d+1}{2}$. Then $\Delta_K(\rho E)$ has positive Lebesgue
measure for almost every $\rho \in SO(d)$.

iii) Suppose that $E$ is radial in the sense that $E=\{r \omega:
\omega \in S^{d-1}; r \in E_0\}$, where $E_0 \subset [0,1]$.
Suppose that the Hausdorff dimension of $E$ is greater than
$\frac{d+1}{2}$. Suppose that $K$ is any symmetric bounded convex
set. Then $\Delta_K(E)$ has positive Lebesgue measure.
\endproclaim

\head Discrete theorems and Continuous $\rightarrow$ Discrete
conversion mechanism \endhead

\vskip.125in

\definition{Definition} We say that $S \subset {\Bbb R}^d$
is well-distributed if there exists a $C>0$ such that every cube
of side-length $C$ contains at least one element of $S$.
\enddefinition

\definition{Definition} We say that a set $S \subset {\Bbb R}^d$
is separated if there exists a constant $c>0$ such that $|a-a'|
\ge c$ for every $a,a' \in A$. \enddefinition

\definition{Definition} Let $K$ be a bounded symmetric convex set
and let $0<\alpha_K \leq d$. We say that the $(K,\alpha_K)$
Falconer conjecture holds if for every compact $E \subset {\Bbb
R}^d$ of dimension greater than $\alpha_K$, $\Delta_K(E)$ has
positive Lebesgue measure. \enddefinition

The essence of the conversion mechanism is captured by the
following result and its proof.

\proclaim{Theorem 0.4} Let $S$ be a well-distributed and separated
subset of ${\Bbb R}^d$. Let $S_q=S \cap {[0,q]}^d$. Suppose that
$(K, \alpha_K)$ Falconer conjecture holds. Then there exists a
constant $c>0$ such that $\# \Delta_K(S_q) \ge
cq^{\frac{d}{\alpha_k}}$.
\endproclaim

In view of Theorem 0.4, every result stated above in Theorem 0.1
and Corollaries 0.2--0.3 has a discrete analog. Moreover, we have
the following application.

\proclaim{Corollary 0.5} Let $S$ be a well-distributed and
separated subset of ${\Bbb R}^d$.

i) Let $K$ be a bounded symmetric convex set. Suppose that
$$ \left| \int_{\partial K} e^{-2 \pi i x \cdot \xi}
d\sigma_K(x) \right| \leq C \gamma({|\xi|}^{-1}), \tag0.4$$ where
$\gamma$ is a convex increasing function with $\gamma(0)=0$ and
$\sigma_K$ is the Lebesgue measure on $\partial K$. Then
$\Delta_K(S)$ is not separated.

ii) Let $K$ be any bounded symmetric convex set. Then
$\Delta_K(\rho S)$ is not separated for almost every $\rho \in
SO(d)$.

\endproclaim

This complements the following result, proved in \cite{Io{\L}a2002} for
$d=2$ and in \cite{Kol2003} for $d\geq 3$.

\proclaim{Theorem 0.6} Let $S$ be well-distributed subset of
${\Bbb R}^d$, and let $\Delta_{K,N}(S)= \Delta_K(S)\cap[0,N]$.

(i) Assume that $d=2$ and $\underline{\lim}_{N\to\infty}\#
\Delta_{K,N}(S)\cdot N^{-3/2}=0$. Then $K$ is a polygon (possibly
with infinitely many sides).
If moreover $\# \Delta_{K,N}(S)\leq CN^{1+\alpha}$ for
some $0<\alpha<1/2$, then the number of sides of $K$ whose length
is greater than $\delta$ is bounded by $C'\delta^{-2\alpha}$.

(iii) Let $d\geq 2$. If $\# \Delta_{K,N}(S)\leq CN$ (in particular,
this holds if $\Delta_K(S)$ is separated), then $K$ is a polytope with
finitely many faces.
\endproclaim

\head Stationary phase tools \endhead

\vskip.125in

In the proofs of our results, we shall make use of the following
estimates on the Fourier transform of the surface carried measure,
which we collect in a single theorem.

\proclaim{Theorem 0.7} Let $K$ be a bounded convex set in ${\Bbb
R}^d$, and let $\sigma_K$ denote the Lebesgue measure on $\partial
K$.

i) Suppose that $\partial K$ is smooth and has everywhere
non-vanishing Gaussian curvature. Then (see e.g. \cite{Herz62})
$$ |\widehat{\sigma}_K(\xi)| \leq C{|\xi|}^{-\frac{d-1}{2}}.
\tag0.5$$

ii) Suppose that $d=2$ and $\partial K$ has everywhere
non-vanishing curvature in the sense of part i) of Corollary 0.3.
Then $(0.5)$ holds without any additional smoothness assumptions.
See, e.g. \cite{BRT98}.

iii) (See \cite{BHI02}) Without any additional assumptions,
$$ {\left( \int_{S^{d-1}} {|\widehat{\sigma}_K (R \omega)|}^2
d\omega \right)}^{\frac{1}{2}} \leq CR^{-\frac{d-1}{2}}, $$
$$ {\left( \int_{S^{d-1}} {|\widehat{\chi}_K (R \omega)|}^2
d\omega \right)}^{\frac{1}{2}} \leq CR^{-\frac{d+1}{2}}. \tag0.6$$

iv) (See \cite{BCT97}) Suppose that $K$ is a polyhedron. Then
$$ \int_{S^{d-1}} |\widehat{\sigma}_K(R \omega)|
d\omega \leq C \log^{d-1}(R) R^{-(d-1)}. \tag0.7$$

\endproclaim

\vskip.125in

\head Proof of Theorem 0.1, Corollary 0.2, and Corollary 0.3
\endhead

\vskip.125in

Let $A_{R, \delta}=\{x \in {\Bbb R}^d: R \leq {||x||}_K \leq
R+\delta\}$ denote an annulus of radius $R$ and width $\delta$,
$\delta\ll R$.

Since the Hausdorff dimension of $E$ is greater than $d-\gamma$,
there is a non-zero Borel measure $\mu$ on $E$ such that the
following energy integral is finite:
$$ \int {|\xi|}^{-\gamma} {|\widehat{\mu}(\xi)|}^2 d\xi<\infty.
\tag1.1$$ For the existence of such a measure, see, for example,
\cite{Falconer85}.

Cover $\Delta_K(E)$ by intervals $\{[R_i,R_i+\delta_i]\}$. It
follows that
$$ 0<(\mu \times \mu)(E \times E) \leq \sum (\mu \times
\mu)\{(x,y): R_i \leq {||x-y||}_K \leq R_i+\delta_i\}$$ $$ \leq C
\sum_i \int \chi_{A_{R_i, \delta_i}}(x-y) d\mu(x)d\mu(y) \leq C
\sum_i \int \widehat{\chi}_{A_{R_i, \delta_i}}(\xi)
{|\widehat{\mu}(\xi)|}^2 d\xi$$ $$ \leq C' \sum_i \delta_i \int
{|\xi|}^{-\gamma} {|\widehat{\mu}(\xi)|}^2 d\xi, \tag1.2$$ where
the last line follows by the decay assumption and the definition
of Lebesgue measure on a hyper-surface. By $(1.1)$, the last
expression in $(1.2)$ is bounded by $C'' \sum_i \delta_i$.

The first part of Theorem 0.1 follows instantly by definition of
measure $0$. The second part follows the same way with point-wise
decay replaced by average decay. The third part follows by an
identical argument.

\vskip.125in

\subhead Proof of Corollary 0.2 \endsubhead The second part of
Corollary 0.2 follows from the first part of Theorem 0.1 and the
first part of Theorem 0.7.  We now prove the first part. We shall
need the following classical result. See, for example,
\cite{BRT98} for a simple proof.

\proclaim{Lemma 1.1} Let $K\subset{\Bbb R}^2$ be a convex body.
Let $\omega=(\cos(\theta), \sin(\theta))$. As before, let
$S_{\theta}=\sup_{x \in K} x \cdot \omega$, and denote by
$l(\theta,\epsilon)$ the length of the chord
$C(\theta,\varepsilon)=\left\{x \in K: x \cdot \omega=S_{\theta}-
\varepsilon \right\}$. Let $\sigma_K$ denote the Lebesgue measure
on $\partial K$. Then, for a constant $C$ independent of
smoothness and curvature, we have

$$ |\widehat{\chi}_K(t\omega)|\leqslant \frac Ct \left(l
\left(\theta,\frac 1{2t}\right)+l\left(-\theta,\frac
1{2t}\right)\right). \tag1.3 $$ \endproclaim

\proclaim{Lemma 1.2} Let $K$ be a convex bounded symmetric set in
the plane. Assume, in addition, that $\partial K$ has
non-vanishing curvature in the sense of Corollary 0.2 (i). Let
$A_{R, \delta}=\{x: R \leq {||x||}_K \leq R+\delta\}$. Then for
$R,|\xi|>1$, $\delta\ll 1$ we have
$$ \left|\widehat{\chi}_{A_{R, \delta}}(\xi) \right| \leq
CR^{\frac{1}{2}} {|\xi|}^{-\frac{1}{2}} \min \{ |\xi|^{-1}, \delta
\}, \tag1.4$$ where $C$ is a constant that depends only on $K$.
\endproclaim

We shall prove Lemma 1.2 in a moment. We first complete the proof
of the first part of Corollary 0.2.

Fix $\beta$ so that $0\leq \beta\leq 1$ and
$\alpha>\frac{1}{2}+\beta$, where $\alpha$ is the Hausdorff
dimension of $E$.  Then there is a non-zero Borel measure $\mu$
supported on $E$ such that the following energy integral is
finite:
$$ \int {|\xi|}^{-\frac{3}{2}+\beta} {|\widehat{\mu}(\xi)|}^2
d\xi<\infty \tag1.5$$ (cf. the proof of Theorem 0.1).

Let $K$ be a convex bounded symmetric planar set satisfying the
assumptions of Corollary 0.3 (i), and let $A_{R,\delta}$ be as
above. We have by Lemma 1.2,
$$ \int \int \chi_{A_{R, \delta}}(x-y) d\mu(x)d\mu(y)=\int
\widehat{\chi}_{A_{R, \delta}}(\xi){|\widehat{\mu}(\xi)|}^2 d\xi$$
$$ \leq CR^{\frac{1}{2}} \left( \int_{|\xi|>\delta^{-1}}
{|\xi|}^{-\frac{3}{2}} {|\widehat{\mu}(\xi)|}^2 d\xi+ \delta
\int_{|\xi| \leq \delta^{-1}} {|\xi|}^{-\frac{1}{2}}
{|\widehat{\mu}(\xi)|}^2 d\xi \right)$$
$$ \leq CR^{\frac{1}{2}} \left(
\delta^\beta \int_{|\xi|>\delta^{-1}} {|\xi|}^{-\frac{3}{2}+\beta}
{|\widehat{\mu}(\xi)|}^2 d\xi
+\delta\cdot\delta^{\beta-1}\int_{|\xi| \leq \delta^{-1}}
{|\xi|}^{-\frac{3}{2}+\beta} {|\widehat{\mu}(\xi)|}^2 d\xi
\right)$$
$$ \leq CR^{\frac{1}{2}} \delta^{\beta}. \tag1.6$$
It follows that
$$ (\mu \times \mu)\{(x,y): R \leq {||x-y||}_K \leq R+\delta\}
\leq CR^{\frac{1}{2}} \delta^{\beta}. \tag1.7$$

Cover $\Delta_K(E)$ by intervals $\{[R_i,R_i+\delta_i]\}$.
Suppose, without loss of generality, that $R_i \leq 10$. It
follows that
$$ 0<(\mu \times \mu)(E \times E) \leq \sum (\mu \times
\mu)\{(x,y): R_i \leq {||x-y||}_K \leq R_i+\delta_i\} \leq C
\sum_i \delta_i^{\beta}. \tag1.8$$

This shows that the $\beta$-dimensional Hausdorff measure of
$\Delta_K(E)$ is non-zero for any
$\beta<\min(1,\alpha-\frac{1}{2})$, so that the Hausdorff
dimension of $\Delta_K(E)$ is at least
$\min(1,\alpha-\frac{1}{2})$.  If furthermore $\alpha>
\frac{3}{2}$, we may take $\beta=1$ and thus deduce that
$\Delta_K(E)$ has positive Lebesgue measure.

We now prove Lemma 1.2. For a fixed $\xi$, let
$$ F(s)=s^2 \widehat{\chi}_K(s \xi). \tag1.9$$

We have
$$ \widehat{\chi}_{A_{R, \delta}}(\xi)=F(R+\delta)-F(R). \tag1.10$$
By the mean value theorem,
$$ |F(R+\delta)-F(R)| \leq \delta \sup_{s \in (R, R+\delta)}
|F'(s)|. \tag1.11$$

Now,
$$ F'(s)=2s \widehat{\chi}_K(s \xi)+s^2 \int_K e^{-2\pi is\xi\cdot x}
(-2\pi i\xi\cdot x)dx =I+II. \tag1.12$$

By Lemma 1.1 and the non-vanishing curvature assumption,
$$ I \leq C s{|s\xi|}^{-\frac{3}{2}}
\leq C {|\xi|}^{-\frac{3}{2}}. \tag1.13$$

On the other hand, following word for word the proof of Lemma 1.1
given in \cite{BRT98} we obtain that
$$ II \leq Cs^2|\xi|{|s\xi|}^{-\frac{3}{2}}
 \leq C R^{\frac{1}{2}}{|\xi|}^{-\frac{1}{2}}. \tag1.14$$

This proves the second estimate in $(1.4)$.  The first estimate
follows from the inequality
$$ |F(R+\delta)-F(R)| \leq |F(R+\delta)|+|F(R)| \leq
CR^2{|R\xi|}^{-\frac{3}{2}} \leq
CR^{\frac{1}{2}}{|\xi|}^{-\frac{3}{2}}, \tag1.15$$ where we again
used Lemma 1.1 and the non-vanishing curvature assumption.

\vskip.125in

\subhead Proof of Corollary 0.3 \endsubhead Part i) follows from
part ii) of Theorem 0.1 and part iv) of Theorem 0.8. The sharpness
result can be obtained as follows. Let $0<s \leq d$. Let $q_1,
q_2, \dots, q_i \dots$ be a sequence of positive integers such
that $q_{i+1} \ge q_i^i$. Let $E_i=\{x \in {\Bbb R}^d: 0 \leq x_j
\leq 1, |x_j-p_j/q_i| \le q_i^{-\frac{d}{s}} \ \text{for some
integers} \ p_j, \ j=1,2\}$. It is not hard to see (see e.g.
\cite{Falconer85}, Chapter 8, or \cite{Wolff02}) that the
Hausdorff dimension of $E=\cap_{i=1}^{\infty} E_i$ is $s$. Also,
$\Delta(E) \subset \bigcap_{i=1}^{\infty} \Delta(E_i)$.

\vskip.125in

Let $K$ be a bounded symmetric convex set in ${\Bbb R}^d$, and let
$P_i=\{p=(p_1,p_2, \dots, p_d): 0 \leq p_j \leq q_i\}$. Then
$\{\|p-p'\|_K:\ p,p'\in P_i\}\subset \{\|p\|_K: p \in P_i\}$, by
translational invariance.  The cardinality of the latter set can
be estimated trivially by $\#P_i\leq (q_i+1)^d$.  We conclude that
$\Delta(E_i)$ is contained in at most $C{(q_i+1)}^d$ intervals of
length bounded by $C' q_i^{-\frac{d}{s}}$. It follows that the Hausdorff
dimension of $\Delta(E)$ is at most $s$. Thus if $s<1$, $\Delta(E)$
has Lebesgue measure $0$.

Part ii) follows from part ii) of Theorem 0.1 and part iii) of
Theorem 0.8.

Part iii) of Corollary 0.3 requires a bit of work. Let $\mu$ be a
probability measure supported on $E$ with the following energy
integral finite:
$$\int {|\xi|}^{-\frac{d-1}{2}} {|\widehat{\mu}(\xi)|}^2 d\xi <\infty.$$
Averaging over rotations if necessary, we may assume that $\mu$ is
rotation-invariant.  Thus $\widehat{\mu}$ is also
rotation-invariant. Let $F(s)=|\widehat{\mu}(s \omega)|$, $\omega
\in S^{d-1}$, and define $A_{R,\delta}$ as in Lemma 1.2. We first
claim that
$$\int_{S^{d-1}}|\widehat{\chi}_{A_{R,\delta}}(r\omega)|d\omega
<C\delta r^{-\frac{d-1}{2}}. \tag1.16$$ Indeed, by $(0.9)$ and
Cauchy-Schwarz we have
$$\int_{S^{d-1}}|\widehat{\chi}_{K}(r\omega)|d\omega
<C r^{-\frac{d+1}{2}}. $$ Now the claim follows by the same
calculation as in the proof of Lemma 1.2, with the above
inequality substituted for Lemma 1.1; the term analogous to $II$
in the proof of Lemma 1.2 is estimated by following the proof of
$(0.9)$ in \cite{BHI02} and using Cauchy-Schwarz again.

We now have
$$ \int \int \chi_{A_{R, \delta}}(x-y) d\mu(x)d\mu(y)=\int
\widehat{\chi}_{A_{R, \delta}}(\xi){|\widehat{\mu}(\xi)|}^2 d\xi$$
$$=\int \int_{S^{d-1}} \widehat{\chi}_{A_{R, \delta}}(r \omega)
{|\widehat{\mu}(r \omega)|}^2 d\omega r^{d-1} dr$$
$$ \leq \int \int |\widehat{\chi}_{A_{R, \delta}}(r
\omega)| d\omega F^2(r)r^{d-1} dr$$
$$ \leq C \delta \int r^{-\frac{d-1}{2}} F^2(r)r^{d-1}dr$$
$$=C \delta \int {|\xi|}^{-\frac{d-1}{2}}
{|\widehat{\mu}(\xi)|}^2 d\xi \leq C'\delta \tag1.17$$ if
$\alpha>\frac{d+1}{2}$ where $\alpha$ is the Hausdorff dimension
of $E$. The fourth line in $(1.17)$ follows by $(1.16)$. The rest
of the proof is exactly like that of Theorem 0.1.

\vskip.125in

\head Proof of Theorem 0.4 and Corollary 0.5. \endhead

\vskip.125in

Let $p=(p_1, \dots, p_d)$. Let $q_i$ denote a sequence of positive
integers such that $q_{i+1} \ge q_i^i$. This sequence will be
specified more precisely below. Let $0<s \leq d$. Let
$$ E_i=\{x \in {[0,1]}^d: |x_j-p_j/q_i| \leq q_i^{-\frac{d}{s}} \
\text{for some} \ \ p \in S\}. \tag2.1$$

Let $E=\bigcap_i E_i$. A standard calculation (see e.g.
\cite{Falconer85}, Ch. 8) shows that the Hausdorff dimension of
$E$ is $s$. Also, $\Delta(E) \subset \bigcap_i \Delta(E_i)$.

Let $S_q$ be defined as in the statement of Theorem 0.2. Suppose
that $\# \Delta_K(S_{q_i}) \leq Cq_i^{\beta}$ for a sequence $q_i$
going to infinity. By refining this sequence we can make sure that
it satisfies the growth condition above. It follows that
$\Delta(E_i)$ can be covered by at most $C q_i^{\beta}$ intervals
of length $C' q_i^{-\frac{d}{s}}$. It follows that the
Hausdorff dimension of $\Delta_K(E)$ is at most
$\frac{s\beta}{d}$.

On the other hand, by assumption, $(K, \alpha_K)$ Falconer
conjecture holds. Letting $s=\alpha_K$, we see that
$\frac{\alpha_K\beta}{d}\geq 1$, so that $\beta \ge
\frac{d}{\alpha_k}$ and we are done.

\vskip.125in

To prove Corollary 0.5 observe that the proof of Theorem 0.4 shows
that if $(0.4)$ holds, then $\frac{\Delta_K(S \cap
{[0,R]}^d)}{R^d} \to \infty$ as $R \to \infty$. This shows that
that $\Delta_K(S)$ cannot be separated.

\newpage

\head References \endhead

\vskip.125in

\ref \key BCT97 \by L. Brandolini, L. Colzani, and G. Travaglini
\paper Average decay of Fourier transforms and integer points in
polyhedra \jour Ark. f. Mat. \vol 35 \yr 1997 \pages 253-275
\endref

\ref \key BIT02 \by L. Brandolini, A. Iosevich, and G. Travaglini
\paper Spherical L1- averages of the Fourier transform of the
characteristic function of a convex set and the geometry of the
Gauss map \yr 2002 \jour Trans. of the AMS (to appear) \endref

\ref \key BHI02 \by L. Brandolini, S. Hofmann, and A. Iosevich
\paper Sharp rate of average decay of the Fourier transform of a
bounded set \jour (submitted) \yr 2002 \endref

\ref \key Bourgain94 \by J. Bourgain \paper Hausdorff dimension
and distance sets \jour Israel J. Math. \vol 87 \yr 1994 \pages
193-201 \endref

\ref \key BRT98 \by L. Brandolini, M. Rigoli, and G. Travaglini
\paper Average decay of Fourier transforms and geometry of convex
sets \jour Revista Mat. Iber. \yr 1998 \vol 14 \pages 519--560
\endref

\ref \key Falconer85 \by K. J. Falconer \paper The geometry of
fractal sets \jour Cambridge University Press \yr 1985 \endref

\ref \key Falconer86 \by K. J. Falconer \paper On the Hausdorff
dimensions of distance sets \jour Mathematika \vol 32 \pages
206-212 \yr 1986 \endref

\ref \key Herz62 \by C. Herz \paper Fourier transforms related to
convex sets \yr 1962 \jour Ann. of Math. \vol 75 \pages 81-92
\endref

\ref \key Io{\L}a2002 \by A. Iosevich and I. {\L}aba \paper Distance
sets of well-distributed planar sets \jour preprint
\yr 2002 \endref

\ref\key Kol2003 \by M. N. Kolountzakis \paper Distance sets 
corresponding to convex bodies \jour (preprint) \yr 2003 \endref

\ref \key Mattila87 \by P. Mattila \paper Spherical averages of
Fourier transforms of measures with finite energy: dimensions of
intersections and distance sets \jour Mathematika \vol 34 \yr 1987
\pages 207-228 \endref

\ref \key Sj\"olin93 \by P. Sj\"olin \paper Estimates of spherical
averages of Fourier transforms and dimensions of sets \yr 1993
\jour Mathematika \vol 40 \pages 322-330 \endref

\ref \key Wolff99 \by T. Wolff \paper Decay of circular means of
Fourier transforms of measures \jour Int. Math. Res. Not. \yr 1999
\vol No.10 \pages 547-567 \endref

\enddocument